\documentclass[12pt]{article}

\usepackage{amsfonts}
\usepackage{amssymb}
\usepackage{enumerate}
\usepackage{amsmath}
\usepackage[T1]{fontenc}
\usepackage{graphicx}

\newtheorem{ttt}{Theorem}[section]
\newtheorem{llll}[ttt]{Lemma}
\newtheorem{ccc}[ttt]{Claim}
\newtheorem{eee}[ttt]{Example}
\newtheorem{fff}[ttt]{Fact}
\newtheorem{rrr}[ttt]{Remark}
\newtheorem{sss}[ttt]{Statement}
\newtheorem{ddd}[ttt]{Definition}
\newtheorem{qqq}[ttt]{Question}
\newtheorem{cccc}[ttt]{Corollary}
\newtheorem{nnn}[ttt]{Notation}
\newtheorem{ppp}[ttt]{Problem}
\newtheorem{ccccc}[ttt]{Conjecture}

\newcommand{\beq}{\begin{equation}}

\newcommand{\bt}{\begin{ttt}}
\newcommand{\bl}{\begin{llll}}
\newcommand{\bc}{\begin{ccc}}
\newcommand{\bex}{\begin{eee}}
\newcommand{\bfa}{\begin{fff}}
\newcommand{\br}{\begin{rrr}\upshape}
\newcommand{\bst}{\begin{sss}}
\newcommand{\bd}{\begin{ddd}\upshape}
\newcommand{\bq}{\begin{qqq}}
\newcommand{\bnn}{\begin{nnn}}
\newcommand{\bpr}{\begin{ppp}}
\newcommand{\bcor}{\begin{cccc}}
\newcommand{\bcon}{\begin{ccccc}}

\newcommand{\eeq}{\end{equation}}
\newcommand{\et}{\end{ttt}}
\newcommand{\el}{\end{llll}}
\newcommand{\ec}{\end{ccc}}
\newcommand{\eex}{\end{eee}}
\newcommand{\efa}{\end{fff}}
\newcommand{\er}{\end{rrr}}
\newcommand{\est}{\end{sss}}
\newcommand{\ed}{\end{ddd}}
\newcommand{\eq}{\end{qqq}}
\newcommand{\ecor}{\end{cccc}}
\newcommand{\econ}{\end{ccccc}}
\newcommand{\enn}{\end{nnn}}
\newcommand{\epr}{\end{ppp}}

\newcommand{\bp}{\noindent\textbf{Proof. }}
\newcommand{\ep}{\hspace{\stretch{1}}$\square$\medskip}

\newcommand{\lab}[1]{\label{#1}}

\newcommand{\NN}{\mathbb{N}}
\newcommand{\ZZ}{\mathbb{Z}}

\newcommand{\RR}{\mathbb{R}}

\newcommand{\PP}{\mathbb{P}}

\newcommand{\al}{\alpha}

\newcommand{\de}{\delta}

\newcommand{\om}{\omega}
\newcommand{\si}{\sigma}

\newcommand{\la}{\lambda}
\renewcommand{\si}{\sigma}
\newcommand{\Si}{\Sigma}

\newcommand{\iL}{\mathcal{L}}

\newcommand{\iF}{\mathcal{F}}

\newcommand{\iN}{\mathcal{N}}

\newcommand{\sm}{\setminus}
\newcommand{\eset}{\emptyset}

\newcommand{\beeq}{\begin{equation}}
\newcommand{\eeeq}{\end{equation}}
\def\su{\subset}

\newcommand{\forces}{\Vdash}
\newcommand{\esm}{\mathfrak{M}}

\newcommand{\height}{\mathrm{ht}}
\renewcommand{\succ}{\mathrm{succ}}
\renewcommand{\phi}{\varphi}

\newcommand{\rt}{\mathrm{root}}

\newcommand{\dom}{\mathrm{dom}}

\numberwithin{equation}{section}

\title{Haar null sets and the consistent reflection of non-meagreness}

\author{
M\'arton Elekes\thanks{Partially supported by the
Hungarian Scientific Foundation grants no.~72655, 61600, 83726 and J\'anos
Bolyai Fellowship.}
\ and 
Juris Stepr\=ans\thanks{Partially supported by NSERC of Canada.}
\\}

\begin{document}

\maketitle

\insert\footins{\footnotesize{MSC codes: Primary 28C10, 03E35  Secondary 03E17,
    22C05, 28A78 }}
\insert\footins{\footnotesize{Key Words: Haar null, Christensen, non-locally
    compact Polish group, packing dimension, Problem FC on Fremlin's list, forcing,
    generic real}}

\begin{abstract}
A subset $X$ of a Polish group $G$ is called \emph{Haar null} if there exists
a Borel set $B \supset X$ and Borel probability measure $\mu$ on $G$ such that
$\mu(gBh)=0$ for every $g,h \in G$.
We prove that there exist a set $X \su \RR$ that is not Lebesgue null and a
Borel probability measure $\mu$ such that $\mu(X + t) = 0$ for every $t \in
\RR$.
This answers a question from David Fremlin's problem list by showing
that one cannot simplify the definition of a Haar null set by leaving out the
Borel set $B$. (The answer was already known assuming the Continuum
Hypothesis.)

This result motivates the following Baire category analogue. It is consistent
with $ZFC$ that there exist an abelian Polish group $G$ and a Cantor
set $C \su G$ such that for every non-meagre set $X \su G$ there exists a $t
\in G$ such that $C \cap (X + t)$ is relatively non-meagre in $C$. This
essentially generalises results of Bartoszy\'nski and Burke-Miller.
\end{abstract}

\section{Introduction}

\subsection{Haar null sets}

Let us first give some motivation for studying Haar null sets in non-locally compact groups.
The following definition is due to Christensen \cite{Ch} (and later
independently to Hunt, Sauer and Yorke \cite{HSY}).

\bd
A subset $X$ of a Polish group $G$ is called \emph{Haar null} if there exists
a Borel set $B \supset X$ and Borel probability measure $\mu$ on $G$ such that
$\mu(gBh)=0$ for every $g,h \in G$.
\ed

(All measures in the paper are assumed to be countably additive.)

The above definition is justified by the following theorem.

\bt[Christensen \cite{Ch}]
A subset of a \emph{locally compact} Polish group is Haar null in the above sense iff it is of Haar measure zero.
\et

There is a huge literature devoted to Haar null sets, e.g. they form a
$\si$-ideal, but Fubini's theorem fails, etc. See for example the work of
Mycielski, Dougherty, Solecki, Matou\v{s}kov\'a, Zaj\'\i \v{c}ek, Duda, Dodos,
Shi, Banakh, Holick\'y, etc. (Note that some
authors use the term \emph{shy} for Haar null, and \emph{prevalent} for
co-Haar null.)

Haar null sets turned out to be very useful in numerous areas of mathematics. 

First, they can serve as a measure counterpart of some Baire category
arguments, even when no actual measure is present. Often the typical behaviour
with respect to this notion of Haar null dramatically differs from the one
with respect to Baire category. For instance,

\bt[Hunt \cite{Hu}]
$\{ f \in C[0,1] : \exists x \ \exists f'(x) \in \RR \}$ is meagre and Haar null.
\et

But:

\bt[Zaj\'\i \v{c}ek \cite{Za}]
$\{ f \in C[0,1] : \exists x \ \exists f'(x) \in [-\infty, \infty] \}$ is meagre but \emph{not} Haar null.
\et

The next two theorems concerning the cycle structure of permutations of the integers also illustrate the striking difference between meagre and Haar null. 

\bt[Folklore]
Comeagre many elements of $S_\infty$ have infinitely many cycles of length $n$ for every $n \in \om$ and no infinite cycles.  
\et

But:

\bt[Dougherty-Mycielski \cite{DM}]
Haar-a.e. element of $S_\infty$ has infinitely many infinite cycles and finitely many finite cycles.
\et

Secondly, Haar null sets show up naturally as exceptional small sets. For example,

\bt[Christensen \cite{Ch2}]
Let $X$ be a separable Banach space and $f : X \to \RR$ a Lipschitz function. Then $f$ is G\^{a}teaux differentiable outside of a Haar null set.
\et

\subsection{Problem FC on Fremlin's list}

After motivating this notion, let us now describe the starting point of the present paper, which is
Problem FC on Fremlin's list\footnote{Originally, Problem FC consisted of
  Problem \ref{p:A} and Problem \ref{p:B}, but after we had solved Problem
  \ref{p:A} and communicated our results to Prof. Fremlin, he has modified the
  problem by mentioning our theorem and erasing the corresponding half of the
  problem.}. The question essentially asks: "But do we need this Borel set $B$
in the definition of Haar null?" The question makes sense for locally compact
groups as well, so it was formulated originally for $\RR$ for simplicity. From
now on we will slightly abuse notation and identify the Borel measure $\mu$
with the outer measure generated by it.

\bpr
\lab{p:A}
Let $X \su \RR$, and let $\la$ denote Lebesgue (outer) measure.
\[
\la(X) = 0 \iff \exists \textrm{ a Borel probability measure } \mu \textrm{ s.t. } \forall t \in \RR \ \mu(X + t) = 0? 
\]
\epr

Note that the left-to-right implication is obvious.
Problem FC also asks if we can find a counterexample $X$ to the right-to-left implication when $\mu = \mu_{Cantor}$ is the usual Cantor measure ("coin tossing measure") on the Cantor set.

\bpr
\lab{p:B}Let $X \su \RR$.
\[
\forall t \in \RR  \ \mu_{Cantor}(X + t) = 0 \implies \la(X) = 0?  
\]
\epr

In fact, Fremlin remarks that the answer to both questions is negative if we assume the Continuum Hypothesis. Let us now prove this for the sake of completeness.

\bc
\lab{c:CH}
Assume the Continuum Hypothesis. Then the answers to Problem \ref{p:A} and Problem \ref{p:B} are in the negative, that is, $\exists X \su \RR$ with  $\la(X) > 0$ such that $\mu_{Cantor}(X + t) = 0$ for every $t \in \RR$.
\ec

\bp
Let $C$ denote the Cantor set. It suffices to construct an $X \su \RR$ with $\la(X) > 0$ such that $C \cap (X+t)$ is countable for every $t \in \RR$.

Let us enumerate the reals as
\[
\{t_\al : \al < \mathfrak{c} \}
\]
and the Borel sets of Lebesgue measure zero as
\[
\{Z_\al : \al < \mathfrak{c} \}.
\]

At stage $\al$ let us pick an
\[
x_\al \in \RR \sm \left( \cup_{\beta < \al} (C - t_\beta) \cup Z_\al \right).
\]

Set
\[
X = \{x_\al : \al < \mathfrak{c}\}.
\]

Then $x_\al \notin Z_\al$ shows that $\la(X) > 0$.
Moreover, for $\al>\beta$, $x_\al + t_\beta \notin C$ implies $C \cap (X + t_\beta) \su \{x_\al : \al \le \beta\}$, hence $C \cap (X + t_\beta)$ is countable.
\ep

\br
In fact we only used $cov(\iN) = cof(\iN)$ (see \cite{BJ} for the
definitions). A more involved argument using ideas similar to the ones in
Section \ref{s:neg} shows that $non(\iN) = \mathfrak{c}$ also suffices.
\er

Therefore the real questions are whether we can find counterexamples in ZFC, that is, without resorting to extra set-theoretic assumptions. 
Our first main goal will be to show in Section \ref{s:neg} (Corollary
\ref{c:FMT}) that Problem \ref{p:A} actually has a negative answer in ZFC.

\bt[First Main Theorem]
\lab{t:FMT}
Problem \ref{p:A} has a negative answer, that is, there exist $X \su \RR$ with 
$\la(X) > 0$ and a Borel probability measure $\mu$ such that $\mu(X + t) = 0$ for every $t \in \RR$. 
\et

Before formulating our second main result, which involves more set theory, let us introduce some terminology.

\bd
Let $X \su \RR$ with $\la(X) > 0$ and $\mu$ be a Borel probability measure. We say that \emph{$\mu$ reflects the positive measure of $X$} if there exists $t \in \RR$ such that $\mu(X + t) > 0$.
\ed

Hence, taking Claim \ref{c:CH} into account, we can reformulate Fremlin's second problem as


\bpr
\lab{p:B2}
Is it consistent that $\mu_{Cantor}$ reflects the positive measure of every $X$?
\epr

This problem is still open, but our second main goal will be to give an affirmative answer to a category analogue in Section \ref{s:pos}.

Now we describe this category analogue in a bit more detail. Theorem \ref{t:FMT} states that there exists a $\mu$ that does not reflect.
Problem \ref{p:B2}, which is still open,  asks if a fixed measure ($\mu_{Cantor}$) consistently reflects. So it is natural to ask the same question about other (fixed) measures.

\bpr
Is it consistent that there exist an \emph{atomless singular} Borel
probability measure $\mu$ such that for every $X \su \RR$ with $\la(X) > 0$
there exists $t \in \RR$ such that $\mu(X + t) > 0$?
\epr

Unfortunately, this is also open. In order to get a better understanding of
the problem, let us consider the following category analogue. Recall that a
set is a Cantor set if it is homeomorphic to the classical middle-thirds
Cantor set.

\bpr
\lab{p:q}
Is it consistent that there exist a Cantor set $K$ such that for every non-meagre set $X \su \RR$ there exists a $t \in \RR$ such that $K \cap (X + t)$ is relatively non-meagre in $K$?
\epr

There are numerous Polish groups that are called "the reals" in set theory. For example, certain results are simpler to prove in $\ZZ_2^\om$ (the Cantor group) than in $\RR$, but usually it is only a technical issue to convert the proofs to $\RR$ (note that the dyadic expansion shows that $[0,1)$ with addition modulo $1$ and $\ZZ_2^\om$ are very similar, the only difference is the presence of carried digits). For technical reasons we will replace $\RR$ with
\[
\RR' = \prod_{m \in \om} \ZZ_{m+3},
\]
where $\ZZ_{m+3} = \{0, 1, \dots, m+2\}$ with addition modulo $m+3$.
Again this group can also be considered as "the reals", since this is essentially an expansion with an "increasing base", as shown by the following map.
\beq
\lab{e:phi}
(n_m)_{m \in \om}  \in \RR' \mapsto \sum_{m \in \om} \frac{n_m}{(m+3)!} \in \RR,
\eeq
which is the analogue of the map $(n_m)_{m \in \om}  \in \ZZ_2^\om \mapsto \sum_{m \in \om} \frac{n_m}{2^m} \in \RR$ that connects the dyadic form of the reals with the usual one.

The second main goal will be to give an affirmative answer to Problem
\ref{p:q} in Section \ref{s:pos} (Theorem \ref{t:SMT}) for this slightly modified underlying group.

\bt[Second Main Theorem]
\lab{t:SMT1}
It is consistent that there exist a Cantor set $C_{EK} \su \RR'$ such that for every non-meagre set $X \su \RR'$ there exists a $t \in \RR'$ such that $C_{EK} \cap (X + t)$ is relatively non-meagre in $C_{EK}$.
\et

Here
\[
C_{EK} = \prod_{m \in \om} \left(\ZZ_{m+3}  \sm \{ m+2 \} \right).
\]
This set was first considered by Erd\H os and Kakutani in \cite{EK}, hence the notation.

\bigskip

The rest of the introduction is devoted to some closely related known results and historical remarks.
If we forget about translates of a fixed Cantor set then the affirmative answer to Problem \ref{p:q} is already known. Interestingly, the following theorem was proved independently in two papers.

\bt[Bartoszy\'nski \cite{Ba}, Burke-Miller \cite{BM}]
It is consistent that for every non-meagre set $X \su \RR$ there exists a
Cantor set $K \su \RR$ such that $K \cap X$ is non-meagre in $K$.
\et

The next theorem still does not use translates, but already finds Cantor sets of some special structure.

\bt[Ciesielski-Shelah \cite{CS}]
For every non-meagre $X \su 2^\NN \times 2^\NN$ there exists a homeomorphism $\varphi : 2^\NN \to 2^\NN$ such that $X \cap graph(\varphi)$ is non-meagre in $graph(\varphi)$.
\et

\br
The Burke-Miller paper also has a certain measure version of their result, and
the Ciesielski-Shelah theorem also has some measure analogue (Ros\l anowski-Shelah \cite{RS}), but these do not seem to say much about Problem \ref{p:q}.
\er

\section{The negative result: Solution to Problem \ref{p:A}}
\lab{s:neg}

The heart of the proof of this result is the following theorem, which is based
on ideas from \cite{DK}. For the definition and basic properties of packing
dimension, denoted by $\dim_p H$, see \cite{Fa} or \cite{Ma}.

\bt
Let $K' \su \RR$ be a Cantor set with $\dim_p K' < 1$ and let $T \su \RR$ be
such that $|T| < \mathfrak{c}$. Then $K' + T$ contains no measurable set of
positive measure.
\et

\bp
Suppose on the contrary that $K' + T$ contains a measurable set $P$ of
positive measure. We may assume that $P$ is compact. By throwing away all
portions (i.e. relatively open nonempty subsets) of measure zero, we may also
assume that every portion of $P$ is of positive measure. In particular, $P$
has no isolated points. The idea of the proof will be to construct a Cantor
set $P' \su P$ such that $P' \cap (K'+r)$ is finite for every $r \in
\RR$. This clearly suffices, since a Cantor set is of cardinality continuum
and hence less than continuum many translates of $K'$ cannot cover $P'$, let
alone $P$.

Let $N$ be a positive integer and let us define $F_N$ to be
the set of $N$-tuples that can be covered by a translate of $K'$, that is,
\[
F_N = \{(x_0,
\dots, x_{N-1}) \in \RR^N : \exists t \in \RR \textrm{ such that } \{x_0,
\dots, x_{N-1} \} \su K' + t \}.
\]
An easy compactness argument shows
that $F_N$ is closed. Reformulating the definition one can easily check that
\[
F_N = (K')^N + \RR(1, \dots, 1),
\]
where $(1, \dots, 1)$ is a vector of $N$ coordinates, and the operations are
Minkowski sum and Minkowski product. It is easy to see that $F_N$ is a
Lipschitz image of $(K')^N \times \RR$, and using that Lipschitz images do not
increase packing dimension as well as $\dim_p (A \times B) \le \dim_p A +
\dim_p B$ and $\dim_p \RR = 1$ we obtain
\[
\dim_p F_N \le N \dim_p K' + 1.
\]
If we choose $N$ large enough, actually if $N > \frac{1}{1 - \dim_p K'}$, then
$N \dim_p K' + 1 < N$, hence 
\[
\dim_p F_N < N.
\]
Let us fix such an $N$.

\bl
\lab{l:shrink}
Let $J_i \su \RR$ $(i<N)$ be closed intervals such that $\mathrm{int}\ \! J_i
\cap P \neq \eset$ $(i<N)$. Then
there are disjoint closed intervals $I_i \su J_i$ $(i<N)$ such that
$\mathrm{int}\ \! I_i \cap P \neq \eset$ $(i<N)$ and
\[
\prod_{i<N} (I_i \cap P) \cap F_N = \eset.
\]
\el

\bp
Since every portion of $P$ is of positive measure, we obtain $\la^N \left(
\prod_{i<N} (\mathrm{int}\ \! J_i \cap P) \right) > 0$, hence $\dim_p
\left(\prod_{i<N} (\mathrm{int}\ \! J_i \cap P) \right) = N > \dim_p
F_N$. Therefore $\left(\prod_{i<N} (\mathrm{int}\ \! J_i \cap P) \right) \sm
F_N \neq \eset$, and, since $F_N$ is closed, $\prod_{i<N} (\mathrm{int}\ \!
J_i \cap P)$ contains a non-empty relatively open set avoiding $F_N$. This
open set contains a basic open set, so there are open intervals $J_i' \su
\mathrm{int}\ \! J_i$ $(i<N)$ intersecting $P$ such that $\prod_{i<N} (J_i'
\cap P) \cap F_N = \eset$.

Finally, since $P$ has no isolated points, it is easy to shrink every $J_i'$
to a closed interval $I_i$ such that they become
disjoint but their interiors still meet $P$. This finishes
the proof of the lemma.
\ep

Now we return to the proof of the theorem. All that remains is to construct
$P'$. We will actually prove
\beq
\lab{e:N}
|P' \cap (K'+r)| < N \textrm{ for every } r \in \RR.
\eeq 
We construct a usual Cantor scheme, where the $k^{th}$ level $\iL_k$ will 
have the following properties for all $k \in \om$.
\begin{enumerate}[(1)]
\item
$\iL_k$ consist of $N^k$ many disjoint closed intervals,
\item
$\forall I \in \iL_{k+1} \exists I'\in \iL_k : I \su I'$,
\item
$\forall I \in \iL_k$ there are $N$ many $I' \in \iL_{k+1}$ with $I' \su I$,
\item
\lab{i:inters}
$\forall I \in \iL_k : \textrm{int}\ \! I \cap P \neq \eset$,
\item
\lab{i:diam}
$\forall I \in \iL_k : \textrm{diam} \ \! I \le \frac{1}{k+1}$,
\item
\lab{i:trans}
 If $I_0, \dots, I_{N-1} \in \iL_k$ are distinct then $\prod_{i<N} (I_i \cap
 P) \cap F_N = \eset$.
\end{enumerate}
(Note that the intervals in (\ref{i:trans}) are not necessarily subsets of the
same $I' \in \iL_{k-1}$.) Assume first that such a Cantor scheme exists, and
define 
\[
P' = \bigcap_{k \in \om} \bigcup \iL_k.
\]
It is easy to see that $P'$ is a Cantor set (\cite{Ke}), while the closedness of
$P$, \eqref{i:inters} and \eqref{i:diam} imply $P' \su P$. Let $x_0, \dots,
x_{N-1}$ be $N$ distinct points in $P'$. Clearly, there is a $k$ and distinct
intervals $I_0, \dots, I_{N-1} \in \iL_k$ such that $x_i \in I_i$
$(i<N)$. Then (\ref{i:trans}) shows that $\{ x_0, \dots,
x_{N-1} \}$ cannot be covered by a translate of $K'$, which proves \eqref{e:N}.

Finally, let us prove by induction that such a Cantor scheme exists. Let
$\iL_0 = \{I\}$, where $I$ is an arbitrary closed interval of length at most $1$
whose interior meets $P$. Assume that $\iL_k$ have already been constructed
with the required properties. Let $\iL_{k+1}'$ be a family of disjoint closed
intervals of length at most $\frac{1}{k+2}$ whose interiors meet $P$ such that
each $I \in \iL_k$ contains $N$ members of $\iL_{k+1}'$. Then recursively
shrinking these intervals by 
applying Lemma \ref{l:shrink} $\binom{N^{k+1}}{N}$ times to all the possible
$N$-tuples of distinct intervals we obtain $\iL_{k+1}$ satisfying all
assumptions. This concludes the proof of the theorem.
\ep

\bt
Let $K \su \RR$ be \emph{a} Cantor set with $\dim_p K < 1/2$. Then there exists $X \su \RR$ with $\la(X) > 0$ such that $|K \cap (X+t)| \le 1$ for every $t \in \RR$.
\et

\bp
As above, let us enumerate the Borel sets of Lebesgue measure zero as $\{Z_\al : \al < \mathfrak{c} \}$.
Since $K - K$ is a Lipschitz image of $K \times K$, we obtain $\dim_p (K - K) < 1$. 
At stage $\al$ let us pick an
\[
x_\al \in \RR \sm \left( \left( (K - K) + \{x_\beta : \beta < \al\} \right)  \cup Z_\al \right).
\]
This is indeed possible by the above theorem applied to $K' = K - K$.
Set 
\[
X = \{x_\al : \al < \mathfrak{c}\}.
\]
Then $x_\al \notin Z_\al$ shows that $\la(X) > 0$.
We still have to check that $|K \cap (X+t)| \le 1$ for every $t \in \RR$. Let
$x_\al, x_\beta \in X$ with $\al > \beta$, and let us assume $x_\al + t,
x_\beta + t \in K$. Then $t \in K - x_\beta$, $x_\al \in K - (K - x_\beta) =
(K-K) + x_\beta$, contradicting the choice of $x_\al$.
\ep

From this we easily obtain our first main theorem (Theorem \ref{t:FMT}) as a corollary.

\bcor[First Main Theorem]
\lab{c:FMT}
Problem \ref{p:A} has a negative answer, that is, there exist $X \su \RR$ with 
$\la(X) > 0$ and a Borel probability measure $\mu$ such that $\mu(X + t) = 0$ for every $t \in \RR$. 
\ecor

\bp
Indeed, let $K$ be any Cantor set with $\dim_p K < 1/2$ (e.g. the "middle-$\al$ Cantor set" is such a set for $\al>1/2$).
Let $\mu$ be any atomless Borel probability measure on $K$. 
Then by the above theorem $\mu$ does not reflect the positive measure of $X$.
\ep

\section{The positive result}
\lab{s:pos}

\subsection{The forcing poset}

The skeleton of the proof of the second main result will be borrowed from Bartoszy\'nski's paper \cite{Ba}.


\bnn
Set $\Si = \bigcup_{l \in \om} \prod_{m<l} \ZZ_{m+3}$,
that is, for $s \in \Si$ the sets $[s] = \{x \in \prod_{m \in \om} \ZZ_{m+3} : s
\su x\}$ form the usual clopen base of $\prod_{m \in \om} \ZZ_{m+3}$.
The symbol $|s|$ will denote the length of the sequence $s \in \Si$, that is,
the cardinality of $\dom(s)$.
\enn
 

Recall that $\forall^\infty$ means `for all but finitely many', and $[\om]^{< \om}$ denotes the set of finite subsets of $\om$.

\bd
Let $s \in \Si$ and $k \in \om$.
Then $S : \om \sm \dom(s) \to [\om]^{< \om}$ is \emph{a finite $k$-slalom above
  $s$}, if 
\begin{enumerate}[(1)]
\item
$\forall i \in \om \sm \dom(s)\ \ |S(i)| \le k$,
\item
$\forall^\infty i \in \om \sm \dom(s)\ \ S(i) = \emptyset$.
\end{enumerate}
\ed

\bd
$\height(S) = \min\{i \in \om \sm \dom(s) : \forall j \ge i \ S(j) =\emptyset \}.$
\ed

\bd
Let $S$ be a finite $k$-slalom above $s$ and $t \in \Si$. Then \emph{$t$
  escapes $S$} if $t \supset s$, $|t| \ge \height(S)$ and $\forall i \in [\dom(s),
  \dom(t))\ \ t(i) \notin S(i)$.
\ed

\bd
Let $s \in \Si$ and $F \su \{t \in \Si : t \supset s\}$. Then \emph{$F$ is
  $k$-fat above $s$}, if for every finite $k$-slalom $S$ above $s$ there
exists $t \in F$ escaping $S$.
\ed

\br
\lab{r:long}
It is easy to see that if $k \ge 1$ and $F$ is $k$-fat above $s$ then for 
every finite $k$-slalom $S$ above $s$ there exist \emph{arbitrarily long} 
$t$'s in $F$ escaping $S$. (Otherwise, just extend $S$ so that $\height(S)$ is
bigger than the length of all $t$'s escaping $S$.)
\er

This immediately yields the following.

\bfa\lab{f:finitefat}
If $F$ is $k$-fat above $s$ and $V$ is a finite set then $F \sm V$ is also
$k$-fat above $s$.
\efa

Recall that $(\Si, \subset)$ is tree, that is, for the purposes of the present
paper, a partially ordered set such that for each $\si \in \Si$ the set
$\{\si' \in \Si : \si' \su \si\}$ is finite. The $n^{th}$ level of a tree is
the set of those points that have exactly $n$ smaller elements. If $\emptyset
\neq T \su \Si$ then $(T, \su)$ is a tree itself.

\bigskip

\bnn
If $t \in T$ then
$\succ_T (t)$ will denote the set of immediate successors of $t$ in $T$. We
simply write $\succ(t)$ when there is no danger of confusion. We say that $T$
has a unique root if it has a unique $\su$-minimal element. In such cases this
root will be denoted by $\rt(T)$. For $t \in T$ let $T[t] = \{s \in T : s
\supset t \}$.
\enn

Now we define our notion of forcing.

\bd
\lab{d:P}
Let $p \in \PP$ iff
\begin{enumerate}[(1)]
\item
$p \su \Si$,
\item
$p$ has a unique root (in particular, $p \neq \eset$),
\item
$\forall t \in p \ \succ_p(t)$ is $1$-fat above $t$.
\item
\lab{k}
$\forall k \in \om \ \forall^\infty t \in p \ \succ_p(t)$ is $k$-fat above $t$.
\end{enumerate}

If $p, p' \in \PP$ then define
\[
p \le_\PP p' \iff p \su p'.
\] 
We will usually simply write $\le$ for $\le_\PP$.
\ed

We will often use the following easy consequence of (\ref{k}).

\bfa
\lab{f:succfat}
If $t \in p \in \PP$ and $k \in \om$ then there exists $s \in p$ such that 
$\succ_p(s)$ is $k$-fat above $s$.
\efa

First we prove that $\PP$ is nontrivial.

\bl
\lab{l:incomparable}
Let $k \ge 1$. If $H$ is $k$-fat above $t$ then $H$ contains a subset
consisting of pairwise incompatible sequences that is $(k-1)$-fat above $t$.
\el

\bp
Let $\{S_n\}_{n \in \om}$ be an enumeration of the $(k-1)$-slaloms above $t$.
It clearly suffices to recursively
pick pairwise incomparable $t_n$'s in $H$ such that $|t_n| > 0$ is strictly increasing and $t_n$ escapes $S_n$.
Suppose $\{ t_m \}_{m<n}$ has already been constructed in such a manner. Then we can form a $k$-slalom by adding
the "last elements of the $t_m$'s" to $S_n$, that is, let
\[
S'_n (i) = 
\begin{cases}
S_n (i) \cup \{ t_m( |t_m|-1 ) \} &  \textrm{if } i = |t_m|-1, \\
S_n (i)                           &  \textrm{otherwise}.\\
\end{cases}
\]

Then $S_n'$ is indeed a $k$-slalom, hence we can choose a $t_n$ escaping
it. By Remark \ref{r:long} we may assume $|t_n| > |t_{n-1}|$. Then the
definition of $S'_n$ shows that $t_n$ is incomparable to $t_m$ for every
$m<n$, and we are done.
\ep

\bl
\lab{l:nonempty}
$\PP \neq \eset$.
\el

\bp
We inductively construct the levels $l_n$ of a tree $p \su \Si$ such that

\begin{enumerate}[(1)]
\item
$|l_0| = 1$,
\item\lab{2}
every $l_n$ consists of pairwise incomparable sequences,
\item\lab{3}
$\forall n \ \forall s \in l_{n+1} \exists t \in l_n \ t \subsetneqq s$, 
\item\lab{4}
$\forall n \ \forall t \in l_n \ l_{n+1} \cap \Si[t]$ is $|t|+1$-fat above
  $t$.
\end{enumerate}
By (\ref{3}) we clearly have $|t| \ge n$ for every $t \in l_n$. Moreover, since $\Si$
is a finitely branching tree,
(\ref{2}) implies that $\forall k \ \forall^\infty t \in l_n \ |t| \ge k$. Using
these two facts and (\ref{4}) it is
easy to see that if such a sequence $\{ l_n \}_n$ exists then $p = \bigcup_n l_n
\in \PP$. 

Let us now check that we can carry out this induction. Suppose that such an
$\{ l_m \}_{m \le n}$ has already been constructed. It is easy to see that
$\Si[t]$ is $|t|+2$-fat above $t$ for every $t \in \Si$. Hence, using Lemma
\ref{l:incomparable}, for every $t \in l_n$ we can pick
$H_t \su \Si[t]$ consisting of pairwise incomparable elements that is 
$|t|+1$-fat above $t$. (We may assume $t \notin H_t$.) Then setting $l_{n+1} =
\bigcup_{t \in l_n} H_t$ completes the proof. 
\ep

The following fact is immediate.

\bfa\lab{f:rest}
If $t \in p \in \PP$ then $p[t] \in \PP$, $p[t] \le p$, and $\rt(p[t]) = t$.
\efa

This easily yields the following three statements.

\bcor
$\PP$ is a separative partial order, and there are incompatible conditions
below every condition. 
\ecor

\bcor\lab{c:fatroot}
$\{ p \in \PP : \succ_p (\rt(p)) \textrm { is } k\textrm{-fat above } \rt(p)  \}$ is dense in $\PP$ for every $k \in \om$.
\ecor

\bcor
$\{p \in \PP : |\rt(p)| \ge k\}$ is dense in $\PP$ for every $k \in \om$. 
\ecor

Next we describe how $\PP$ adds a generic real. The last corollary easily
implies that if $G \su \PP$ is a generic filter then $\dot{x}_G = \bigcup_{p \in G}
\rt(p)$ is a function $\dot{x}_G \in \prod_{m \in \om} \ZZ_{m+3}$. From now on we
denote by $\dot{x}$ a name for this generic real.

\br
Some textbooks require that forcing posets have largest elements, but
our $\PP$ has no such element. One possible answer to this problem is that one
can actually do forcing without largest elements (since we can basically `add
a largest element to $\PP$'), and hence some other textbooks actually avoid largest
elements in the definition of forcing posets. But there is another possible
answer in case of $\PP$; by mimicking the proof of Lemma \ref{l:nonempty} it is
not hard to see that $\PP$ is dense in
\[
\PP_0 = \{ p \su \Si  :\forall k \in \om \ \forall^\infty t \in p \ p[t]
\textrm{  is } k\textrm{-fat above } t \},
\]
which already has a largest element, namely $\Si$. 
\er
The reason why we prefer $\PP$ to the apparently simpler $\PP_0$ is that it
fits our fusion arguments (inductive constructions) better.

\subsection{Properness and preservation of non-meagreness}

It will be necessary, of course, to prove that $\PP$ is proper, but for the
intended iteration we will need a stronger property. Recall that $p' \in \PP$
is $\esm$-generic, if for
every dense open $D \su \PP$ if $D \in \esm$ then $p' \forces `` D \cap \dot{G} \cap
\esm \neq \eset"$. (Here $\dot{G}$ is a name for the generic filter.)

\bd 
A forcing notion  $\PP$ is said to be
\emph{Cohen-preserving} if for every condition $p \in \PP$, 
every countable elementary submodel $\mathfrak{M}$ such that $p,\PP, \le_\PP
 \in \mathfrak{M}$ and every real $c$ that is a Cohen over 
$\mathfrak{M}$, there is an $\mathfrak{M}$-generic condition $p' \leq p$ such
 that  $p' \forces ``c \text{ is a Cohen real over } \mathfrak{M}[\dot{G}]"$. 
\ed

We now spell out this last clause in a bit more detail. For more information
see e.g. \cite{BJ}. Let $\mathcal{F} = \{f: \Si \to \Si \ | \ \forall \si \in
\Si \ f(\si) \supset \si \}$. Then the dense open subsets of $\prod_{m \in \om}
\ZZ_{m+3}$ are precisely the sets of the form $U_f = \bigcup_{\si \in \Si}
[f(\si)]$, where $f$ ranges over $\iF$. Then $p' \forces ``c \text{ is a Cohen real over }
\mathfrak{M}[\dot{G}] "$ means that if $\dot{f}$ is a name for an element of
$\mathcal{F}$ and $\dot{f} \in \mathfrak{M}$ then $p' \forces ``c \in
U_{\dot{f}} \ "$. 

It is not hard to see that the notions of Cohen-preserving and second category
set preserving partial order does not depend on the underlying Polish
space. We will only use that $\prod_{m \in \om} \ZZ_{m+3}$ and $\RR$ are the same in
this respect, which follows e.g. from the fact that we can throw away
countable sets from these spaces so that the remaining sets are
homeomorphic. (See the map in \eqref{e:phi} in the discussion preceding Theorem \ref{t:SMT1}.)

The following results are well-known, see e.g. \cite{BJ}.

\bt\lab{t:shelah} 
Cohen-preserving partial orders are proper and they preserve second category
sets. The countable support iteration of Cohen-preserving partial orders is also
Cohen-preserving. 
\et

\bl\lab{l:q^*}
Let $p^* \in \PP$, $\dot{f}^*$ be a name for an element of $\mathcal{F}$ and
$D^*$ be a dense open subset of $\PP$ such that $p^*, \dot{f}^*, D^* \in \esm$.
Then there exists $q^* \le p^*$ such that $q^* \in \esm \cap D^*$ and $q^*
\forces ``c \in U_{\dot{f}^*} \ "$. 
\el

\bp
Define
\[
V = \bigcup \{ [\si'] : \exists q^* \le p^* \ \exists \si \in \Si \ q^* \forces ``
\dot{f}^*(\si) = \si'\  " \}.
\]
It is not hard to see that $V \su \prod_{m \in \om} \ZZ_{m+3}$ is dense open and $V
\in \esm$. Since $c$ is Cohen over $\esm$, we obtain $c \in V$, so we can find
$q^*, \si$ and $\si'$ such that $c \in [\si']$ and $q^* \forces
``\dot{f}^*(\si) = \si' \ "$. Let us now fix $\si$ and $\si'$, then clearly
$\exists q^* \in D^* \ q^* \le p^* \ q^* \forces ``\dot{f}^*(\si) = \si'
\ "$. Applying elementarity to this last formula we obtain such a $q^* \in
\esm$. Since $q^* \forces ``\dot{f}^*(\si) = \si' \ "$ clearly implies $q^*
\forces ``c \in U_{\dot{f}} \ "$, the proof is complete.
\ep

Now we are ready to prove the main result of the section. The proof will essentially be an inductive construction of a condition. Unlike in the proof of Lemma \ref{l:nonempty}, we will not build the tree `level-by level', but we will perform a kind of `back-and forth' fusion instead. The only place where this more complicated fusion is essential is Lemma \ref{l:non-meager}, but we decided to use this method here in a simpler situation as well to make the reading of Lemma \ref{l:non-meager} easier.

\bl\lab{l:Cohen}
$\PP$ is Cohen-preserving (and hence proper, as well).
\el

\bp
Let $\{ \dot{f}_n \}_{n \in \om}$ enumerate the names for elements of $\mathcal{F}$ that are in $\esm$, and
let $\{ D_n \}_{n \in \om}$ enumerate the dense open subsets of $\PP$ that are in $\esm$.

For $n \in \om$ we will inductively define
\begin{enumerate}[(i)]
\item
$s_n \in \Si$,
\item\lab{t_n}
$q_n \in \PP$,
\item
$t_n \in \Si$,
\item
$p_n \in \PP$,
\end{enumerate}

\noindent
such that for every $m \le n$ the following hold:

\begin{enumerate}[(1)]
\item
\lab{i:in1}
$t_m \in p_n$, 
\item
\lab{i:succ1}
$\succ_{p_m}(t_m) \sm \{s_0, \dots, s_n\} \su \succ_{p_n}(t_m)$,
\item
\lab{i:m+11}
$\succ_{p_m}(t_m)$ is $(m+1)$-fat above $t_m$,
\item
\lab{i:decr1}
$p \ge p_0$ and $p_m \ge p_n$,
\item
\lab{i:q11}
$q_m = p_m [t_m]$,
\item
\lab{i:q21}
$q_n \in D_m$,
\item
\lab{i:q31}
$q_n \forces ``c \in U_{\dot{f_m}} \ "$.
\end{enumerate}

We will make sure that every stage of the induction will be carried out in $\esm$, and we will tacitly assume that all object we pick at the stages are in $\esm$. (The whole induction will of course not be in $\esm$.) 

Let us start with $n=0$. Put $s_0 = \eset$. By Lemma \ref{l:q^*}, Corollary \ref{c:fatroot} and Fact \ref{f:rest} there exists $q_0 \le p$ such that $q_0 \in D_0$, $q_0 \forces ``c \in U_{\dot{f_0}} \ "$, and if $t_0 = \rt(q_0)$ then $\succ_{q_0}(t_0)$ is $1$-fat above $t_0$. Setting $p_0 = q_0$ finishes the $0^{th}$ step. It is not hard to check that the inductive assumptions are satisfied. 

Let us now assume that $s_m, q_m, t_m$ and $p_m$ have already been defined for $m \le n$ satisfying the inductive assumptions. For every $m$ let 
$\{S_m^k\}_{k \in \om}$ be an enumeration of the set of $(m+1)$-slaloms above $t_m$. To start the $n+1^{st}$ step, first we need to pick a $t_m$ for some $m \le n$. We make sure by some simple bookkeeping that during the course
of the induction each $t_m$ will be picked infinitely many times, and when we
visit the node $t_m$ for the $k^{th}$ time then we take care of $S_m^k$ (we
construct a $t_{n+1}$ above $t_m$ escaping $S_m^k$).

So let us assume that we are at the $n+1^{st}$ step and we pick $t_m$ for the $k^{th}$ time. Inductive assumption (\ref{i:m+11}) yields that $\succ_{p_m} (t_m)$ is $(m+1)$-fat above $t_m$, hence so is $\succ_{p_m} (t_m) \sm \{s_0, \dots, s_n\}$ by Fact \ref{f:finitefat}.
Thus we can fix an $s_{n+1} \in \succ_{p_m} (t_m)\sm \{s_0, \dots, s_n\}$ escaping the $(m+1)$-slalom $S_m^k$. By (\ref{i:succ1}) we have $s_{n+1} \in p_n$ as well. Applying
Lemma \ref{l:q^*} $n+2$ times we obtain a $q_{n+1} \le p_n [s_{n+1}]$ such that
$q_{n+1} \in \bigcap_{m \le n+1} \ D_m$ and $q_{n+1} \forces ``c \in \bigcap_{m \le n+1} \ U_{\dot{f_m}} \ "$. By Corollary \ref{c:fatroot} and Fact \ref{f:rest} we may assume that if $t_{n+1} = \rt(q_{n+1})$ then $\succ_{q_{n+1}}(t_{n+1})$ is $n+2$-fat above $t_{n+1}$. Setting $p_{n+1} = (p_n \sm p_n[s_{n+1}] ) \cup q_{n+1}$ finishes the $n+1^{st}$ step.

Now we check that the inductive assumptions are satisfied. Items (\ref{i:in1}) and (\ref{i:succ1}) follow from the structure of the fusion. Namely, it is not hard to see that at the $n+1^{st}$ step we only modify $p_n$ in the `cone' $p_n [s_{n+1}]$, and this cone does not contain the earlier $t_m$'s, moreover, an element of $\succ_{p_m} (t_m)$ only `disappears' from $p_n$ when it is picked as an $s_{n+1}$.
Items (\ref{i:decr1}), (\ref{i:q11}), (\ref{i:q21}) and (\ref{i:q31}) are straightforward from the construction, and (\ref{i:m+11}) follows from (\ref{i:q11}).

Let us now define $p' = \{ t_m \}_{m \in \om}$. It is easy to see that $p' \in
\PP$, since we made sure by the bookkeeping that $\succ_{p'}(t_m)$ is $m+1$-fat above $t_m$ for every
$m$. Combining (\ref{i:in1}) and (\ref{i:decr1}) we obtain 
\beq
\lab{e:p'1}
p' \le p_n
\eeq 
for every $n$, and also that $p' \le p$.

All that remains to be shown is that $p'$ is $\esm$-generic and $p' \forces ``c \text{ is a Cohen real over } \mathfrak{M}[\dot{G}] "$.

First we check that $p'$ is $\esm$-generic. Let $n_0 \in \om$ be fixed, then we have
to show that $p' \forces `` D_{n_0} \cap \dot{G} \cap \esm \neq \eset "$. Let $p''
\le p'$ be arbitrary, it suffices to find a $p''' \le p''$ such that $p'''
\forces `` D_{n_0} \cap \dot{G} \cap \esm \neq \eset "$. Since every condition is
infinite, there exists $n \ge n_0$ such that $t_n \in p''$. By (\ref{i:q21}) we have 
$q_n \in D_{n_0}$ and also $q_n \in \esm$, so $q_n \forces ``
D_{n_0} \cap \dot{G} \cap \esm \neq \eset "$. Thus, \eqref{e:p'1} and 
(\ref{i:q11}) imply $p''[t_n] \le
p'[t_n] \le p_n[t_n] = q_n$, so we are done by choosing $p'''= p''[t_n]$.

A similar argument shows that $p' \forces ``c \text{ is a Cohen real
  over } \mathfrak{M}[\dot{G}] "$. Indeed, for every $n_0$ and $p'' \le p'$
there exists $n \ge n_0$ such that $t_n \in p''$. Then $q_n \forces ``c \in
U_{\dot{f}_{n_0}}"$, $q_n \in \esm$, and $p''' = p''[t_n] \le
p'[t_n] \le p_n [t_n] = q_n$, so we are done.
This finishes the proof of the lemma.
\ep

\subsection{The Main Lemma}

Our main lemma will describe how a single step in the iterated forcing
construction works.

\bl
Let $X \in V$ such that $X \su \prod_{m \in \om} \ZZ_{m+3}$ and $X$ is non-meagre in every
non-empty open subset of $\prod_{m \in \om} \ZZ_{m+3}$. Then $\PP
\forces `` X \cap (\dot{C}_{EK} - \dot{x}) \textrm{ is non-meagre in }
\dot{C}_{EK} - \dot{x} "$.
\el

\bp
It is not hard to see that if a set $H \su C_{EK} - x$ is meagre
in $C_{EK} - x$ then there exists a decreasing sequence of open sets $U_n \su \prod_{m \in
  \om} \ZZ_{m+3}$ such that $U_n \su B(C_{EK} - x, \frac{1}{n+1})$, $U_n \cap
(C_{EK} - x)$ is dense in $C_{EK} - x$ and $\bigcap_n U_n \cap H = \emptyset$.
(Here $B(A, \varepsilon)$ denotes the $\varepsilon$-neighbourhood of the set $A$.)

Hence let us assume that there exist $p \in \PP$ and a name $ \{ \dot{U}_n \}_{n \in \om}
$ for a decreasing sequence of open subsets of $\prod_{m \in \om} \ZZ_{m+3}$ such that $p
\forces `` \dot{U}_n \su B(\dot{C}_{EK} - \dot{x}, \frac{1}{n+1})$, $\dot{U}_n
\cap (\dot{C}_{EK} - \dot{x})$ is dense in $\dot{C}_{EK} - \dot{x}$ and
$\bigcap_n \dot{U}_n \cap X = \emptyset "$.

Define
\[
R_{p,  \{ \dot{U}_n \}_{n \in \om} } = \{r \in \prod_{m \in \om} \ZZ_{m+3} : \exists p' \le p,
p' \forces `` r \in \bigcap_n \dot{U}_n"\}.
\]

For the definitions and basic facts concerning analytic sets, sets with the
property of Baire, etc, we refer the reader to \cite{Ke}.

In the next subsection we will prove that we can assume (by replacing $p$ with
a stronger condition, if necessary) that $R_{p, \{ \dot{U}_n \}_{n \in \om} }$ is
analytic. Therefore it possesses
the property of Baire. Moreover, in Subsection \ref{s:non-meager} we will prove that $R_{p,
\{ \dot{U}_n \}_{n \in \om} }$ is non-meagre. Let us now
accept these statements for the moment. Then $R_{p, \{ \dot{U}_n \}_{n \in \om} }$ is
co-meagre in a non-empty open set, thus $X \cap R_{p, \{ \dot{U}_n \}_{n \in \om} } \neq
\eset$, so we can fix an $r \in X \cap R_{p, \{ \dot{U}_n \}_{n \in \om} }$. But then $r
\in X$ and $p' \forces `` r \in \bigcap_n \dot{U}_n"$ for some $p' \le p$, thus
$p' \forces ``\bigcap_n \dot{U}_n \cap X \neq \emptyset "$. On the other hand,
$p' \le p$ implies that $p' \forces ``\bigcap_n \dot{U}_n \cap X = \emptyset "$,
which is a contradiction.
\ep

\bcor
\lab{c:main}
Let $X \in V$, $X \su \prod_{m \in \om} \ZZ_{m+3}$, $X$ is non-meagre. Then $\PP
\forces `` \exists t \in \prod_{m \in \om} \ZZ_{m+3} \textrm{ such that } X \cap
(\dot{C}_{EK} + t) \textrm{ is non-meagre in } \dot{C}_{EK} + t "$.
\ecor

\bp
Let $Q$ be the analogue of the rationals, that is, $Q = \{ q \in \prod_{m \in
  \om} \ZZ_{m+3} : \forall^\infty m \ q(m) = 0\}$. Then $X + Q$ is non-meagre in
every non-empty 
open subset of $\prod_{m \in \om} \ZZ_{m+3}$, hence $\PP \forces ``(X + Q) \cap
(\dot{C}_{EK} - \dot{x}) \textrm{ is non-meagre in } \dot{C}_{EK}  - \dot{x}"$
by the previous lemma. But $Q$ is
countable and $X + Q = \bigcup_{q \in Q} X + q$, so $\PP \forces `` \exists
q \in Q \textrm{ such that } (X + q) \cap (\dot{C}_{EK} - \dot{x})
\textrm{ is non-meagre in } \dot{C}_{EK} - \dot{x}"$. But then $\PP \forces
``X \cap (\dot{C}_{EK} - \dot{x} - q) \textrm{ is non-meagre in }
\dot{C}_{EK} - \dot{x} - q"$ and we are done.
\ep

We still have to prove the two statements concerning $R_{p, \{ \dot{U}_n \}_{n
    \in \om} }$.

\subsubsection{Analyticity of $R_{p, \{ \dot{U}_n \}_{n \in \om} }$}

As above, let $p \in \PP$ and $ \{ \dot{U}_n \}_{n \in \om} $ be a name for a decreasing sequence
of open sets of $\prod_{m \in \om} \ZZ_{m+3}$ such that $p \forces `` \dot{U}_n \su
B(\dot{C}_{EK} - \dot{x}, \frac{1}{n+1})$, $\dot{U}_n \cap (\dot{C}_{EK} -
\dot{x})$ is dense in $\dot{C}_{EK} - \dot{x}$ and $\bigcap_n \dot{U}_n \cap X
= \emptyset "$. Recall that

\[
R_{p,  \{ \dot{U}_n \}_{n \in \om} } = \{r \in \prod_{m \in \om} \ZZ_{m+3} : \exists p' \le p,
p' \forces `` r \in \bigcap_n \dot{U}_n"\}.
\]

\bl
\lab{l:analytic}
There exists a $q \le p$ such that $R_{q,  \{ \dot{U}_n \}_{n \in \om} }$ is analytic.
\el

We will split the proof into several steps.

\bd
Let $p \in \PP$. A set $B \su p$ is called a \emph{barrier} if it intersects
every infinite branch of $p$. It is \emph{open}, if $s,t \in p$, $s \su t$, $s \in B$
imply $t \in B$. 
\ed

\bd
Let $p$ and $\{ \dot{U}_n \}_{n \in \om}$ be as above. Then $p$ is \emph{nice with
respect to $\{ \dot{U}_n \}_{n \in \om}$} if for every $s \in \Si$ and every $n \in
\om$
\[
B_{s,n} = \{t \in p : p[t] \forces ``[s] \su \dot{U}_n" \textrm{ or } p[t]
``\forces [s] \not\su \dot{U}_n" \}
\]
is a barrier. (It is clearly open.)
\ed

Note that if $B \su p$ is an open barrier and $q \le p$ then $q \cap B \neq
\eset$.

\bl
Assume that $p$ is nice with respect to $\{ \dot{U}_n \}_{n \in \om}$. Let $q^*
\le p$, $n^* \in \om$, $r \in \prod_{m \in \om} \ZZ_{m+3}$ and $q^* \forces
``r \in \dot{U}_{n^*}"$. Then there are $t^* \in q^*$ and $k^* \in \om$ such that
$p[t^*] \forces ``[r|{k^*}] \su \dot{U}_{n^*}"$.
\el

\bp
There are $q' \le q^*$ and $k^* \in \om$ such that $q' \forces ``[r|{k^*}] \su
\dot{U}_{n^*}"$. Since $B_{r|{k^*}, n^*}$ is an open barrier in $p$ and $q' \le p$,
we obtain that there is a $t^* \in q' \cap B_{r|{k^*}, n^*}$. Then $p[t^*] \forces
``[r|{k^*}] \not\su \dot{U}_{n^*}"$ is impossible, since then $q'[t^*] \le p[t^*], q'$
would force two contradicting statements. Hence $p[t^*] \forces ``[r|{k^*}] \su
\dot{U}_{n^*}"$ by the definition of $B_{r|{k^*}, n^*}$, and we are done.
\ep

\bl
\lab{l:B_n}
Assume that $p$ is nice with
respect to $\{ \dot{U}_n \}_{n \in \om}$. Let $r \in \prod_{m \in \om} \ZZ_{m+3}$ be
arbitrary. Then $r \in R_{p,  \{ \dot{U}_n \}_{n \in \om} }$ iff there exists a sequence
of sets $B_n \su p$ and a sequence of functions $\phi_n : B_n \to \om$ such
that for every $n \in \om$
\begin{enumerate}[(1)]
\item
$|B_0| = 1$,
\item
$B_n$ consists of pairwise incomparable sequences,
\item
$\forall t \in B_{n+1} \ \exists s \in B_n \ s \subsetneqq t$,
\item
$\forall s \in B_n \ \phi_n(s) > n$,
\item
$\forall k \in \om \ \forall^\infty s \in B_n \ \phi_n(s) \ge k$,
\item
$\forall s \in B_n \ \{ t \in B_{n+1} : t \supset s \}$ is
  $\phi_n(s)$-fat above $s$,
\item
\lab{7}
$\forall t \in B_n \ \exists k \in \om \ p[t] \forces ``[r|k] \su
  \dot{U}_n"$.
\end{enumerate}
\el

\bp
Let us first suppose $r \in R_{p,  \{ \dot{U}_n \}_{n \in \om} }$, that is, $p' \forces ``r \in \bigcap_n \dot{U}_n"$ for some $p' \le p$. We construct 
$\{ B_n \}_{n \in \om}$ and $\{\phi_n\}_{n \in \om}$ by induction on $n$. We will make the induction hypothesis that
\[
\forall s \in B_n \ \succ_{p'} (s) \textrm{ is } \phi_n(s)\textrm{-fat above }
s.
\]

Applying the previous lemma with $q^* = p'$ and $n^* = 0$ we obtain $t_0 \in
p'$ and $k_0 \in \om$ such that $p[t_0] \forces ``[r|k_0] \su
\dot{U}_0"$. Put $B_0 = \{t_0\}$. By Fact \ref{f:succfat} we may
assume (by extending $t_0$ in $p'$ if necessary) that $\succ_{p'}(t_0)$ is
$1$-fat above $t_0$. Define $\phi_0 (t_0) = 1$. Then all requirements imposed
on $B_0$ and $\phi_0$ are satisfied.

Suppose that $B_n$ and $\phi_n$ satisfying all requirements have
already been constructed. For each $t \in \bigcup_{s \in B_n}
\succ_{p'} (s)$ let us apply the previous lemma with $q^* = p'[t]$ and $n^* =
n+1$. Thus we obtain $u_t \in p'$, $u_t \supset t$ and $k_t \in \om$ such that
$p[u_t] \forces ``[r|k_t] \su \dot{U}_{n+1}"$. Let $\psi : \bigcup_{s \in B_n}
\succ_{p'} (s) \to \om$ be an arbitrary function satisfying
\[
\forall t \in \bigcup_{s \in B_n} \succ_{p'} (s) \ \  \psi(t) > n + 1
\]
and
\[
\forall k \in \om \ \forall^\infty t \in \bigcup_{s \in B_n} \succ_{p'} (s)
\ \ \psi(t) \ge k.
\]
Using Fact \ref{f:succfat} we may assume that $u_t$ is $\psi(t)$-fat above
$u_t$ for every $t$. (Note that extending $u_t$ in $p'$ if necessary does not
harm any of the requirements.) Then putting $B_{n+1} = \{ u_t : t \in \bigcup_{s \in
B_n} \succ_{p'} (s) \}$ and $\phi_{n+1} (u_t) = \psi(t)$ finishes the proof
of this direction.

In order to prove the other direction, let us assume that $\{ B_n \}_{n \in
\om}$ and  $\{\phi_n\}_{n \in \om}$ satisfy all requirements. Set $p' = \bigcup_n
B_n$. It is easy to see that $p' \in \PP$ and $p' \le p$. Therefore it
suffices to check that for every $n_0 \in \om$ we have $p' \forces ``r \in
\dot{U}_{n_0}"$. Otherwise, there exists a $p'' \le p'$ such that
\beq
\lab{e:1}
p'' \forces ``r \notin \dot{U}_{n_0}".
\eeq
It is easy to see from the construction of $p'$ that $B_{n_0}' = \{s \in p' : \exists t \in B_{n_0} \ \ t \su
s \}$ is an open barrier in $p'$, hence $p'' \cap B_{n_0}' \neq \eset$. Thus there are
$t \in B_{n_0}$ and $s \supset t$ with $s \in p''$. Then $p''[s] \le p[t],
p''$, and $p[t] \forces ``r \in \dot{U}_{n_0}"$ by (\ref{7}), which is a
contradiction by (\ref{e:1}).
\ep

\bl
Let $p$ and $\{ \dot{U}_n \}_{n \in \om}$ be as above. Then there exists $q
\le p$ that is nice with respect to $\{ \dot{U}_n \}_{n \in \om}$.
\el

\bp
The proof will be similar to the previous inductive construction, so we will
omit some details. Let $\{ (s_i, n_i) \}_{i \in \om}$ be an enumeration of
$\Si \times \om$. We inductively define a sequence $B_i \su p$ and for every
$t \in B_i$ a condition $p_t \le p$ with $\rt(p_t) = t$ as follows.

The fact `$q \forces \phi$ or $q \forces \neg \phi$' will be abbreviated as $q ||
\phi$.

There exists $p' \le p$ such that $p' || ``[s_0] \su \dot{U}_{n_0}"$. We may assume
that $\succ_{p'} (\rt(p'))$ is $1$-fat above $\rt(p')$. Let $t_0
= \rt(p')$ and put $B_0 = \{ t_0 \}$, $p_{t_0} = p'$. 

Now assume that $B_i$ and $p_t$ for every $t \in B_i$ have already been
constructed. For every $t \in B_i$ and every $u \in \succ_{p_t} (t)$ find $p'
\le p_t[u]$ such that $p' || ``[s_{i+1}] \su \dot{U}_{n_{i+1}}"$ and $\succ_{p'}
( \rt(p') )$ is `sufficiently fat above $\rt(p')$'. Then let $B_{i+1}$ be the
set of these $\rt(p')$'s and let $p_{\rt(p')} = p'$. This finishes the general
step of the induction. Note that $p_t || ``[s_i] \su \dot{U}_{n_i}"$ for every $i
\in \om$ and every $t \in B_i$.

Then it is not hard to see that $q = \bigcup_i B_i \in \PP$ and $q \le
p$. It easily follows from the construction that every $B_i$ is a barrier in
$q$. Let us now check that $q$ is nice with respect to $\{ \dot{U}_n \}_{n \in
\om}$. So let us fix $i \in \om$, and it suffices to show that $q[t] ||
``[s_i] \su \dot{U}_{n_i}"$ for every $t \in B_i$. But this is clear, since $p_t
|| ``[s_i] \su \dot{U}_{n_i}"$ and $q[t] \le p_t$.
\ep

Now we are ready to prove Lemma \ref{l:analytic}.

\bp
By the previous lemma we may assume that $p$ is nice with respect to $\{
\dot{U}_n \}_{n \in \om}$. Define
\[
B_{p, \{ \dot{U}_n \}_{n \in \om}} =
\]
\[
\left\{ \left( r, \{B_n\}_{n \in \om}, \{ \phi_n
\}_{n \in \om} \right) : r, \{B_n\}_{n \in \om} \textrm{ and } \{ \phi_n
\}_{n \in \om} \textrm{ are as in  Lemma \ref{l:B_n}} \right\}
\]
\[
\su  \prod_{m \in \om} \ZZ_{m+3} \times (2^p)^\om \times (2^{p
  \times \om} )^\om,
\]
where $\ZZ_{m+3}, \om$ and $p$ are all equipped with the discrete topology, hence
this huge ambient product space is compact metric, therefore Polish. It suffices to prove that $B_{p, \{
  \dot{U}_n \}_{n \in \om}}$ is a Borel set in this product space, since then
$R_{p, \{ \dot{U}_n \}_{n \in \om}}$, which is the projection of $B_{p, \{
  \dot{U}_n \}_{n \in \om}}$ on the first coordinate, is analytic. We mostly
leave this standard but very lengthy computation to the reader, and only deal
with the most interesting clause, that is, Lemma \ref{l:B_n} (\ref{7}). 

The conditions `$\forall n \in \om$', `$\forall t \in p$', `$t \in B_n$',
`$\exists k \in \om$' are clearly Borel, so it suffices to check that for
fixed $n, t$ and $k$
\[
V_{n, t, k} = \{ r \in \prod_{m \in \om} \ZZ_{m+3} : p[t] \forces `` [r|k] \su \dot{U}_n"   \}
\] 
is Borel. But clearly
\[
V_{n, t, k} =  \bigcup \{ [s] : s \in \Si, |s| = k, p[t] \forces `` [s] \su
\dot{U}_n" \},
\]
so it is actually a union of basic clopen sets, hence open. 
\ep

\subsubsection{Non-meagreness of $R_{p, \{ \dot{U}_n \}_{n \in \om} }$}
\lab{s:non-meager}

Yet again, let $p \in \PP$ and $ \{ \dot{U}_n \}_{n \in \om} $ be a name for a decreasing sequence
of open sets of $\prod_{m \in \om} \ZZ_{m+3}$ such that $p \forces `` \dot{U}_n \su
B(\dot{C}_{EK} - \dot{x}, \frac{1}{n+1})$, $\dot{U}_n \cap (\dot{C}_{EK} -
\dot{x})$ is dense in $\dot{C}_{EK} - \dot{x}$ and $\bigcap_n \dot{U}_n \cap X
= \emptyset "$. Recall that

\[
R_{p,  \{ \dot{U}_n \}_{n \in \om} } = \{r \in \prod_{m \in \om} \ZZ_{m+3} : \exists q \le p,
q \forces `` r \in \bigcap_n \dot{U}_n"\}.
\]

\bl
\lab{l:q^**}
Let $q^* \le p$, $s^* \in \Si$ and $l^* \in \om$ such that $q^* \forces
``[s^*] \cap (\dot{C}_{EK} - \dot{x}) \neq \eset"$. Then there exist $q^{**}
\le q^*$ and $s^{**} \in \Si, s^{**} \supset s^*$ such that $q^{**}\forces
``[s^{**}] \cap (\dot{C}_{EK} - \dot{x}) \neq \eset \textrm{ and } [s^{**}]
\su \dot{U}_{l^*}"$.
\el

\bp
Set $\dot{C} = \dot{C}_{EK} - \dot{x}$. Since $p \forces ``\dot{U}_{l^*} \cap
\dot{C} \textrm{ is dense open in } \dot{C}"$ and $q^* \forces ``[s^*] \cap
\dot{C} \textrm{ is non-empty open in } \dot{C}"$ we obtain that
$q^* \forces ``[s^*] \cap \dot{U}_{l^*} \cap \dot{C}  \neq \eset"$. 
Since $[s^*] \cap \dot{U}_{l^*}$ is a name for an open subset of $\prod_{m \in
  \om} \ZZ_{m+3}$, and if an open set meets a set then it contains a basic open
set meeting the same set, we infer that  
$q^* \forces ``\exists \dot{s}^{**} \in \Si, \ [\dot{s}^{**}] \cap \dot{C}
\neq \eset, \ [\dot{s}^{**}] \su [s^*] \cap \dot{U}_{l^*}"$.
Hence there exist $q^{**} \le q^*$ and $s^{**} \in \Si, s^{**} \supset
s^*$ such that $q^{**}\forces ``[s^{**}] \cap \dot{C} \neq \eset \textrm{ and
} [s^{**}] \su \dot{U}_{l^*}"$.
\ep

\bd
We will write $s||t$ to denote that $s(i) + t(i) \neq i+2$ for every $i <
\min(|s|, |t|)$.

For $p \in \PP$ and $t \in \Si$ the symbol $p||t$ will abbreviate that $s||t$ for every $s \in p$.
\ed

The following fact can be easily checked by a standard argument using that $p \forces \varphi \iff \forall p'\le p \ \exists p''  \le p'\ p'' \forces \varphi$. The details are left to the reader.

\bfa
\lab{f:||}
Let $p \in \PP$ and $t \in \Si$. Then the following are equivalent.
\begin{enumerate}
\item
$p||t$,
\item
$p \forces ``[t] \cap (\dot{C}_{EK} - \dot{x}) \neq \eset"$.
\end{enumerate}
\efa

We will also need one more lemma. The proof, which is left to the reader again, follows easily from the definition of fatness.

\bl
\lab{l:prune}
Let $k>0$, $s_0 \in \Si$, $F \su \{s \in \Si : s \supset s_0\}$ be
$k$-fat above $s_0$, and also let $\si \in \Si$ such that $s_0 || \si$. Then $\{t \in F : t||\si \}$
is $k-1$-fat above $s_0$.
\el



We are now ready to prove what we are aiming at.

\bl
\lab{l:non-meager}
$R_{p,  \{ \dot{U}_n \}_{n \in \om} }$ is non-meagre.
\el

\bp
We have to show that $R_{p,  \{ \dot{U}_n \}_{n \in \om} }$ intersects every
dense $G_\de$ set, so it suffices to prove that $R_{p,  \{ \dot{U}_n \}_{n \in
    \om} } \cap \bigcap_n V_n = \eset$ for every
sequence $\{ V_n \}_{n \in \om}$ of dense open subsets of $\prod_{m \in \om}
\ZZ_{m+3}$. 
The proof will work as follows. On the one hand, we will inductively define a
strictly increasing sequence $\{r_n\}_{n \in \om}$ of elements of $\Si$ 
such that $[r_n] \su V_n$,
which will of course imply that if we set $r = \bigcup_n r_n$ then $r \in
\bigcap_n V_n$. On the other hand, we will also simultaneously carry
out a fusion argument similar to the one in the proof of Lemma \ref{l:Cohen}
and obtain a $p' \le p$ such that $p'\forces `` r \in \bigcap_n
\dot{U}_n"$. This will show $r \in R_{p,  \{ \dot{U}_n \}_{n \in \om} } \cap
\bigcap_n V_n $, which will complete the proof. 

Let us now start the fusion. The main differences between this argument and the
one in Lemma \ref{l:Cohen} will be that we will use Lemma \ref{l:q^**} instead
of Lemma \ref{l:q^*}, and we will be building the $r_n$'s as well.

For $n \in \om$ we will inductively define
\begin{enumerate}[(i)]
\item
$s_n \in \Si$,
\item
$q_n \in \PP$,
\item
$t_n \in \Si$,
\item
$r_n \in \Si$,
\item
$p_n \in \PP$,
\end{enumerate}

\noindent
such that for every $m \le n$ the following hold:

\begin{enumerate}[(1)]
\item
\lab{i:||}
$r_n || t_m$,
\item
\lab{i:long}
$|r_n| \ge |t_m|$, 
\item
\lab{i:in}
$t_m \in p_n$, 
\item
\lab{i:succ}
$\succ_{p_m}(t_m) \sm \{s_0, \dots, s_n\} \su \succ_{p_n}(t_m)$,
\item
\lab{i:m+2}
$\succ_{p_m}(t_m)$ is $(m+2)$-fat above $t_m$,
\item
\lab{i:decr}
$p \ge p_0$ and $p_m \ge p_n$,
\item
\lab{i:q1}
$q_m = p_m [t_m]$,
\item
\lab{i:q2}
$q_m \forces ``[r_m] \su \dot{U}_m"$.
\end{enumerate}

Let us start with $n=0$. Put $s_0 = \eset$. Applying Lemma \ref{l:q^**} to
$p$, $\eset$ and $0$, then using Corollary \ref{c:fatroot} and Fact
\ref{f:rest}, we obtain that there exists $q_0 \le
p$ and $r_0' \in \Si$ such that $q_0\forces ``[r_0'] \cap (\dot{C}_{EK} -
\dot{x}) \neq \eset
\textrm{ and } [r_0'] \su \dot{U}_0"$, and if $t_0 = \rt(q_0)$ then
$\succ_{q_0}(t_0)$ is $2$-fat above $t_0$. Then $r_0' || t_0$ by Fact
\ref{f:||}. Hence we can clearly find an $r_0'' \supset r_0'$ such that
$|r_0''| \ge |t_0|$ and $r_0'' || t_0$. Finally, since $V_0$ is dense open, we
can extend $r_0''$ further to obtain an $r_0 \supset r_0''$ with $[r_0] \su
V_0$. Setting $p_0 = q_0$ finishes the $0^{th}$ step. It is not hard to check
that the inductive assumptions are satisfied. (Note that (\ref{i:||}) and
(\ref{i:long}) follow from $|r_0''| \ge |t_0|$, $r_0'' || t_0$, and $r_0
\supset r_0''$.)

Let us now assume that $s_m, q_m, t_m, r_m$ and $p_m$ have already been
defined for 
$m \le n$ satisfying the inductive assumptions. For every $m$ let $\{S_m^k\}_{k \in
  \om}$ be an enumeration of the set of $(m+1)$-slaloms above $t_m$.
  To start the $n+1^{st}$ step, first we need to pick a $t_m$ for some
$m \le n$. We make sure by some simple bookkeeping that during the course
of the induction each $t_m$ will be picked infinitely many times, and when we
visit the node $t_m$ for the $k^{th}$ time then we take care of $S_m^k$ (we
construct a $t_{n+1}$ above $t_m$ escaping $S_m^k$).

So let us assume that we are at the $n+1^{st}$ step and we pick $t_m$ for the
$k^{th}$ time. Inductive assumption (\ref{i:m+2}) yields that $\succ_{p_m}
(t_m)$ is $(m+2)$-fat above $t_m$, hence so is $\succ_{p_m} (t_m) \sm \{s_0,
\dots, s_n\}$ by Fact \ref{f:finitefat}. Therefore $\{t \in \succ_{p_m} (t_m)
\sm \{s_0, \dots, s_n\} : t || r_n \}$ is $(m+1)$-fat above $t_m$ by \eqref{i:||} and Lemma
\ref{l:prune}. Thus, using Remark
\ref{r:long} as well, we can fix a $s_{n+1} 
\in \succ_{p_m} (t_m)\sm \{s_0, \dots, s_n\}$ escaping the $(m+1)$-slalom
$S_m^k$ such that $s_{n+1} || r_n$, and $|s_{n+1}| \ge |r_n|$. By (\ref{i:succ})
we also have $s_{n+1} \in p_n$. As $s_{n+1} || r_n$ and $|s_{n+1}| \ge
|r_n|$, we obtain $p_n [s_{n+1}] || r_n$. Hence Lemma \ref{l:q^**} applied to
$p_n [s_{n+1}]$, $r_n$ and $n+1$, then Corollary \ref{c:fatroot} and Fact
\ref{f:rest} yield a $q_{n+1} \le p_n[s_{n+1}]$ and a
$r_{n+1}' \supset r_n$ such that $q_{n+1} || r_{n+1}'$ and
$q_{n+1} \forces ``[r_{n+1}'] \su \dot{U}_{n+1}"$,
and if $t_{n+1} = \rt(q_{n+1})$ then
\beq
\lab{e:q_{n+1}}
\succ_{q_{n+1}}(t_{n+1}) \textrm{ is }  n+3\textrm{-fat above } t_{n+1}.
\eeq 
Then $r_{n+1}' || t_{n+1}$, hence we can clearly find an $r_{n+1}'' \supset
r_{n+1}'$ such that $|r_{n+1}''| \ge |t_{n+1}|$ and $r_{n+1}'' ||
t_{n+1}$. Finally, since $V_{n+1}$ is dense open, we
can extend $r_{n+1}''$ further to obtain an $r_{n+1} \supset r_{n+1}''$ with
$[r_{n+1}] \su V_{n+1}$. Setting $p_{n+1} = (p_n \sm p_n [s_{n+1}]) \cup
q_{n+1}$ finishes the $n+1^{st}$ step. 

Now we check that the inductive assumptions are satisfied. For (\ref{i:||}) and (\ref{i:long})
it suffices to check that $r_{n+1} || t_{n+1}$ and $|r_{n+1}| \ge |t_{n+1}|$, which is analogous to the case $n=0$ above. 
Items (\ref{i:in}) and (\ref{i:succ}) follow from the structure of the fusion, as already described in Lemma \ref{l:Cohen}.
Namely, at the $n+1^{st}$ step we only modify $p_n$ in the `cone' $p_n [s_{n+1}]$, and this cone does not contain the earlier $t_m$'s, moreover, an element of $\succ_{p_m} (t_m)$ only `disappears' from $p_n$ when it is picked as an $s_{n+1}$. Items (\ref{i:decr}), (\ref{i:q1}), and (\ref{i:q2}) are straightforward from the construction, and (\ref{i:m+2}) follows from (\ref{i:q1}) and \eqref{e:q_{n+1}}.

Let us now define $p' = \{ t_m \}_{m \in \om}$. It is easy to see that $p' \in
\PP$, since $\succ_{p'}(t_m)$ is $m+1$-fat above $t_m$ for every
$m$. Combining (\ref{i:in}) and (\ref{i:decr}) we obtain 
\beq
\lab{e:p'}
p' \le p_n
\eeq 
for every $n$, and also that $p' \le p$.

What remains to be shown is that $p' \forces ``r \in \bigcap_n \dot{U}_n"$,
that is, $p' \forces ``r \in \dot{U}_{n_0}"$ for every fixed $n_0$. Let $p''
\le p'$ be arbitrary, then it suffices to find a $p''' \le p''$ such that $p'''
\forces ``r \in \dot{U}_{n_0}"$.  As every condition is infinite, there exists
$n \ge n_0$ such that $t_n \in p''$. Then \eqref{e:p'},
(\ref{i:q1}), and (\ref{i:q2}) imply $p''[t_n] \le p'[t_n] \le
p_n[t_n] = q_n \forces ``[r_n] \su \dot{U}_n"$. Therefore, since the $U_n$'s
are decreasing and $r \in [r_n]$, we obtain 
$p''[t_n] \forces ``r \in \dot{U}_{n_0}"$. Thus $p''' = p''[t_n]$
works, and this finishes the proof of the lemma.
\ep

\subsection{Putting the proof together}

\bt[Second Main Theorem]
\lab{t:SMT}
It is consistent with $ZFC$ that for every non-meagre set $X \su \prod_{m \in
\om} \ZZ_{m+3}$ there is some $t \in \prod_{m \in \om} \ZZ_{m+3}$ such that $X \cap
(C_{EK} + t)$ is non-meagre in $C_{EK} + t$.
\et

\bp
Iterate $\PP$ of length $\omega_2$ with countable support over a model $V$
of the Continuum Hypothesis to obtain $V_\alpha$ for $\alpha \leq \omega_2$.
If $ V_{\omega_2} \models X \su \prod_{m \in \om} \ZZ_{m+3} \textrm{ is
non-meagre}$, then by an easy reflection argument there is an $\al < \om_2$
such that $V_\alpha \models X \cap V_\alpha \text{ is non-meagre}$ (see the
analogous \cite[Lemma 12.]{Ba}). Applying Corollary \ref{c:main} yields that in
$V_{\alpha+1}$ there is some $t_{\al + 1} \in \prod_{m \in \om} \ZZ_{m+3}$ such that
$V_{\al+1} \models (X \cap V_\al) \cap (C_{EK} + t_{\al + 1}) \textrm{ is
non-meagre in } C_{EK} + t_{\al + 1}$.  Then Theorem~\ref{t:shelah} implies
that
$V_{\omega_2} \models (X\cap V_{\alpha}) \cap (C_{EK} + t_{\al + 1}) \textrm{
is non-meagre in } C_{EK} + t_{\al + 1}$.  Hence for the larger set we also
obtain $V_{\omega_2} \models X \cap (C_{EK} + t_{\al + 1}) \textrm{ is
non-meagre in } C_{EK} + t_{\al + 1}$, which completes the proof.
\ep

\section{Open problems}

In this final section we collect the open questions.

\bpr
Let $X \su \RR$.
\[
\forall t \in \RR  \ \mu_{Cantor}(X + t) = 0 \implies \la(X) = 0?  
\]
\epr

\bpr
Is it consistent that there exist an \emph{atomless singular} Borel
probability measure $\mu$ such that for every $X \su \RR$ with $\la(X) > 0$
there exists $t \in \RR$ such that $\mu(X + t) > 0$?
\epr

\bpr
Is $\PP$ of Definition \ref{d:P} forcing equivalent to the Miller forcing?
\epr

\bpr
Does Theorem \ref{t:SMT} hold for $\RR$ instead of $\prod_{m \in \om} \ZZ_{m+3}$?
\epr

\bigskip

\noindent
\textsc{R\'enyi Alfr\'ed Institute, Re\'altanoda u. 13-15. Budapest
1053, Hungary\\
and\\
Institute of Mathematics, E\"otv\"os Lor\'and University, 
P\'azm\'any P\'eter s. 1/c, Budapest 1117, Hungary}

\textit{Email address}: \verb+emarci@renyi.hu+

{\tt www.renyi.hu/\hbox{$\sim$}emarci}

\bigskip

\noindent
\textsc{Department of Mathematics, York University,
Toronto, Ontario M3J 1P3, Canada }

\textit{Email address}: \verb+steprans@yorku.ca+

\end{document}